\documentclass{article}
\usepackage{amsmath,amsfonts,epsfig}
\usepackage{mathptmx}

\usepackage{amsmath,amsfonts,latexsym,amssymb}
\usepackage{mathrsfs}
\usepackage{verbatim}
\usepackage[utf8]{inputenc}
\usepackage[T1]{fontenc}
\usepackage{ae,aecompl}
\usepackage{braket}

\newcommand{\bb}{\mathbb}

\newcommand{\cx}{{\bb C}}

\newcommand{\hthree}{{\bb H}^3}
\newcommand{\htwo}{{\bb H}^2}

\renewcommand{\bold}[1]{\medskip \noindent {\bf \boldmath #1
                        }\nopagebreak[4]}

\newcommand{\chat}{\widehat{\cx}}

\newcommand{\area}{\operatorname{area}}

\newcommand{\inj}{\operatorname{inj}}

\renewcommand{\Re}{\operatorname{Re}}

\newcommand{\Teich}{\operatorname{Teich}}

\newcommand{\tr}{\operatorname{tr}}
\newcommand{\vol}{\operatorname{vol}}

\newtheorem{theorem}{Theorem}[section]
\newtheorem{prop}[theorem]{Proposition}
\newtheorem{lemma}[theorem]{Lemma}
\newtheorem{cor}[theorem]{Corollary}

\newcommand{\cL}{{\cal L}}
\newcommand{\cM}{{\cal M}}

\renewcommand{\hbar}{\bar{{\mathbb H}}^3}

\newcommand{\CC}{\mathbb C}

\newcommand{\R}{\mathbb R}

\newcommand{\Rvol}{V_R}
\newcommand{\Cvol}{V_C}

\newcommand{\psl}{{\rm PSL}(2,\CC)}

\renewcommand{\htwo}{{\mathbb H}^2} 
\def\eproof{$\Box$ \medskip}
\renewcommand{\area}{{\operatorname{{\bf area}}}}
\newcommand{\dome}{{\operatorname{{\bf dome}}}}
\newcommand{\ho}{{\mathfrak{h}}}

\renewcommand\marginpar[1]{} 

\begin{document}

\title{\bf \Large Schwarzian derivatives, projective structures, and the Weil-Petersson gradient flow for
  renormalized volume 
} \author{Martin Bridgeman\thanks{Research 
    supported by NSF grant DMS-1500545} \ and Jeffrey
  Brock\thanks{Research supported by NSF grant DMS-1608759.} \ and
  Kenneth Bromberg\thanks{Research supported by NSF grant
    DMS-1207873.}}

\date{\today}

\maketitle

\begin{abstract}
 To a complex projective structure $\Sigma$ on a surface, Thurston
  associates a locally convex pleated surface. We derive bounds on the
  geometry of both in terms of the norms $\|\phi_\Sigma\|_\infty$ and
  $\|\phi_\Sigma\|_2$ of the quadratic differential $\phi_\Sigma$ of
  $\Sigma$ given by the Schwarzian derivative of the associated
  locally univalent map. We show that these give 
  a unifying approach that generalizes a number of important,
  well known results for convex cocompact hyperbolic structures on
  3-manifolds, including bounds on the Lipschitz constant for the
  nearest-point retraction and the length of the bending
  lamination. We then use these 
  bounds to begin a study of the Weil-Petersson gradient flow of
  renormalized volume on the space $CC(N)$ of convex cocompact
  hyperbolic structures on a compact manifold $N$ with incompressible
  boundary, leading to a proof of the conjecture that the renormalized
  volume has infimum given by one-half the simplicial volume of $DN$,
  the double of $N$.
\end{abstract}

\section{Introduction}

Throughout the work of Bers, Sullivan, and Thurston, the precise
relation between the
conformal boundary of a hyperbolic 3-manifold and its internal geometry has
been a key subtlety. For example, the classical {\em Bers inequality}
bounds the lengths of geodesics in the 3-manifold in terms of their
lengths in the hyperbolic metric on the conformal boundary, and a
related theorem of Sullivan gives uniform bounds on the Teichm\"uller
distance between the conformal boundary and the boundary of the convex
core for 3-manifolds with incompressible boundary. A long history of 
 results of this type, obtained by Canary \cite{CAN91}, Bishop
\cite{Bish}, Epstein-Marden-Markovic \cite{EMM1} and Bridgeman-Canary
\cite{BC05,   BC10}, have made important advances through
delicate arguments.

This paper provides a unifying perspective to these considerations via
the {\em Schwarzian derivative}, which naturally associates a
holomorphic quadratic differential to each component of the conformal
boundary of a hyperbolic 3-manifold. Remarkably, in addition to
shining new light on a number of important results in the literature,
the `Schwarzian' is key to proving a conjectured lower bound on the {\em
  renormalized volume} of hyperbolic 3-manifolds, a notion whose import we
elucidate here.

To begin with, the following initial result illustrates these explicit
connections:
\begin{theorem}\label{summary}
Let $M$ be a hyperbolic 3-manifold, $\partial_c M$ its conformal
boundary and $C(M)$ its convex core. Let $\phi_M$ be the holomorphic
quadratic differential obtained from the Schwarzian derivative of the map
comparing $\partial_c M$ to its Fuchsian uniformization. Then 
\begin{enumerate}
\item the retract map $\partial_c M\to \partial C(M)$ is
  $\sqrt{1+2\|\phi_M\|_\infty}$-Lipschitz, and
\item $L(\lambda_M) \le 4\pi|\chi(\partial_c
  M)|\|\phi_M\|_\infty$ where $L(\lambda_M)$ is the length of the  bending
  lamination $\lambda_M$ of $\partial C(M)$.
\end{enumerate}
\end{theorem} 
Indeed, Theorem~\ref{summary} follows almost directly from a theorem
of G. Anderson bounding Thurston's {\em projective metric} in terms of
the hyperbolic metric where the bound depends on the Schwarzian
derivative. Taking Anderson's result together with the classical
Nehari bound on the Schwarzian, we obtain many well known
results, such as the Lipschitz bounds of Epstein-Marden-Markovic
\cite{EMM1} and Bridgeman-Canary \cite{BC10}, and the length bounds of
Bridgeman-Canary \cite{BC05}, as immediate corollaries.

Working a bit harder, we obtain bounds on $L(\lambda_M)$ in terms of
the $L^2$-norm of the Schwarzian, which we employ to study the
powerful notion of {\em renormalized volume}. Motivated by
considerations from theoretical physics, the notion of 
renormalized volume was first introduced by Graham and Witten (see
\cite{GW}) in the general setting of conformally compact Einstein
manifolds. In the setting of infinite volume, convex cocompact
hyperbolic 3-manifolds, renormalized volume has been seen to be of
particular interest as a more analytically natural proxy for convex core
volume (see for example \cite{TT, TZ}). The approach here follows the
work of Krasnov-Schlenker \cite{KS08} and Schlenker
\cite{schlenker-qfvolume}. Our $L^2$-bounds give the following tight
relationship between the convex core volume $V_C(M)$ and the
renormalized volume $V_R(M)$ of $M$.

\begin{theorem}\label{L2_squeeze}
There is a function $G(t) \sim t^{1/5}$ such that  if $M$ is a convex
cocompact hyperbolic three-manifold with incompressible boundary
then 
$$\Cvol(M) - |\chi(\partial M)|G(||\phi_M||_2) \leq \Rvol(M) \leq \Cvol(M)$$
and $\Rvol(M) = \Cvol(M)$ if and only if the convex core of $M$ has
totally geodesic boundary. 
\end{theorem}

The result reveals the close connection of the renormalized volume to the
volume the convex core, but renormalized volume carries the advantage
that if we fix a hyperbolizable 3-manifold $N$ then $\Rvol$ is a
smooth function on the space $CC(N)$ of all convex cocompact
hyperbolic 3-manifolds homeomorphic to $N$, and Krasnov-Schlenker have
established a formula for its derivative (see \cite{KS08} and Theorem
\ref{Rvolvariation}). It is natural to conjecture that the infimum
$\mathcal V_R(N)$ of $\Rvol$ is the purely topologically defined {\em
  simplicial volume} of $N$. By applying the variational formula of
Krasnov-Schlenker and Theorem \ref{L2_squeeze} to study the
Weil-Petersson gradient flow of $\Rvol$, we establish the conjectured
lower bound:
%
\begin{cor}\label{simplicial_vol}
Let $N$ be a compact hyperbolizable 3-manifold with non-empty
incompressible boundary and without torus boundary components. Then
$\mathcal V_R(N) = \frac12 V_S(DN)$ where $DN$ is the double of $N$
and $V_S(DN)$ is the simplicial volume. The infimum is realized if and
only if $N$ is acylindrical or has the homotopy type of a closed
surface. 
\end{cor}
Corollary~\ref{simplicial_vol} is an analogue of a well-known 
result of Storm on the convex core volume (see \cite{storm2}).  

Partial results in this direction were established
prior to our work. It follows immediately from the Krasnov-Schlenker
variational formula that all critical points of $\Rvol$ occur at $M
\in CC(N)$ where the convex core of $M$ has totally geodesic
boundary. Note that this can only occur when $N$ is acylindrical, in
which case there is a unique such structure in $CC(N)$, or $N$ is
homotopy equivalent to a surface and there is a half dimensional
subspace of $CC(N)$ of Fuchsian structures where the renormalized
volume is zero. In the acylindrical setting, Moroianu (see \cite{SM})
and Palette (see \cite{PAL2}) have independently shown that this
critical point is a local minimum of $\Rvol$. When $N$ is homotopy
equivalent to a closed surface, the ``surface group" case, our result implies that $\mathcal V_R(N)
=0$. Previously, Krasnov-Schlenker (see \cite{KS08}) had proven that renormalized volume has zero infimum when taken over quasifuchsian manifolds with finitely bent convex core boundary. 
In the special case of {\em almost-Fuchsian}
structures this was proven by Ciobotaru-Moroianu (see
\cite{CM}). Finally, when $N$ is acylindrical, Corollary
\ref{simplicial_vol} was proven by Palette \cite{PAL} using very
different methods. In fact, combining our methods with those of
Palette gives a new and technically simpler proof of the Storm theorem
on convex core volume for acylindrical manifolds. Note that prior to
the work here it was not even known that the renormalized volume was
positive.  \medskip

\noindent{\bf Core volume, renormalized volume, and Weil-Petersson
  distance.} In a sequel \cite{BBB}, we study the Weil-Petersson gradient flow
further, supplying a direct proof of renormalized
volume lower bounds in terms of Weil-Petersson distance: 
\begin{theorem}\label{rvol_lower}
Given $\epsilon > 0$ there exists a $c = c(\epsilon,S) > 0 $ so that
if $d_{\rm WP}(X,Y) \geq \epsilon$ then we have
$$\Rvol(Q(X,Y)) \geq c \cdot d_{WP}(X,Y).$$
\end{theorem}
Here, $Q(X,Y)$ denotes the {\em Bers simultaneous uniformization} of $X$ and
$Y$ in $\Teich(S)$, and $d_{\rm WP}(X,Y)$ is their {\em Weil-Petersson
distance}.
Together with the comparison of Theorem~\ref{L2_squeeze}, we obtain
direct proofs of the lower bounds on convex core volume in
\cite{brock-wp} and \cite{Brock:3ms1}.  Previously, these results had
been obtained by building a combinatorial model for the Weil-Petersson
metric (the {\em pants graph}) and showing these combinatorics also
give volume estimates for the relevant convex cores.  The model relies
on delicate combinatorial arguments involving the complex of curves
and its hierarchical structure developed in
\cite{Masur:Minsky:CCI,Masur:Minsky:CCII}, while the renormalized
volume flow produces a natural analytic link between Weil-Petersson
distance and volume.


\paragraph{Outline.} 
We begin with a  discussion of locally univalent maps and complex projective structures. On a
projective structure there are two natural metrics: the hyperbolic
metric, which depends only on the underlying conformal structure, and
Thurston's {\em projective metric}.  By comparing a projective
structure to its Fuchsian uniformization one also obtains a
holomorphic quadratic differential via the Schwarzian derivative. The
main technical tool of the paper is an unpublished theorem of
G. Anderson (Theorem \ref{thurston_metric_bound}) bounding the
projective metric in terms of the hyperbolic metric and a function of
the $L^\infty$-norm of the Schwarzian quadratic
differential. Section 2.1 is devoted to a short, new proof of this
theorem. As with the original, the proof is based on a construction of
Epstein which associates a surface in $\hthree$ to a conformal
metric on the unit disk $\Delta$ and a locally univalent map
$f\colon\Delta\to\chat$.

In Sections 2.2 and 2.3  we review Thurston's parameterization of locally
univalent maps and of projective structures in terms of measured
laminations. In particular Thurston parameterizes projective
structures on a surface by locally convex pleated surfaces. There is
a natural ``retract'' map from the projective structure to the pleated
surface that is 1-Lipschitz from the projective metric to the path
metric on the pleated surface. Using the Schwarzian bound on the
projective metric we obtain a bound on the Lipschitz constant for the
retract map when we take the hyperbolic metric on the domain (Corollary \ref{projective_lipschitz}). The
length of the bending lamination is also controlled by the Schwarzian
as it is a linear function of the area of the projective metric (Theorem \ref{infinity_norm_bending_bound}). 

In Section 2.4 we review the classical bounds of Nehari on the $L^\infty$-norm of the Schwarzian derivative of univalent maps and use the Nehari bounds to bound the Schwarzian when the locally univalent map is a covering map for a domain in $\Omega$. In Sections 2.5 and 2.6 we combine the Nehari bounds to derive Lipschitz bounds on retract maps from domains in $\chat$ to convex hulls in $\hthree$ (Theorem \ref{dome_lipschitz}) and from the conformal boundary of a hyperbolic 3-manifold to the boundary of the convex core (Theorem \ref{convex_core_lipschitz}). We also obtain bounds on the length of the bending lamination of the convex core (Theorems \ref{inj_bending_bound} and \ref{strong_compressible_bound}). 

All of these bounds are based on the $L^\infty$-norm of the Schwarzian. In Section 2.7 we bound the length of the bending lamination in terms of the $L^2$-norm of the Schwarzian. This will be used in our study of renormalized volume.

In Section 3, the last part of the paper, we begin our study of the renormalized
volume of a convex cocompact hyperbolic 3-manifold. After reviewing definitions we improve on bounds, originally due the Schlenker, comparing the renormalized volume to the volume of the convex core. In particular, we show that the difference of the two volumes is bounded by a function of the $L^2$-norm of the Schwarzian of the projective boundary (Theorem \ref{L2_squeeze}). 

We use these bounds to study Weil-Petersson gradient flow of $-V_R$. Along flow lines, the $L^2$-norm of the Schwarzian of the projective boundary will decay to zero. It will follow that the infimum of renormalized volume will agree with the infimum of convex core volume (Theorem \ref{renorm_inf}).

We highlight one other novelty of our approach, a new definition of the {\em $W$-volume}. The
usual definition of $W$-volume involves the integral of the mean
curvature over the boundary of the manifold. We'll see that it can be
reinterpreted as a function of the volume of the submanifold, the area
of the boundary and the area of its associated metric at
infinity. This reinterpretation is valid even when the boundary is not
smooth and clarifies the formula for the $W$-volume of the convex core given in \cite{schlenker-qfvolume}.

The proof of our theorem on the lower bound for renormalized volume is actually quite short. The reader who is solely interested in this result can skip much of the paper as it only depends on the bound on the projective metric (Theorem \ref{thurston_metric_projective}), the bound on the length of the bending lamination in terms of the $L^2$-norm of the Schwarzian (Section 2.7), and Section 3.

{\bf Acknowledgements:} The authors would like to thank Dick Canary for comments and suggestions on the paper. We would also like to thank the referee whose comments helped improve the paper greatly.

\section{Epstein surfaces and projective structures}
\newcommand{\Ep}{\operatorname{Ep}}
Let $f\colon \Delta\to \chat$ be a locally univalent map. Thurston
defined a natural metric on $\Delta$ associated to $f$, the {\em
  Thurston} or {\em projective metric}. Here is the definition: Let $D
\subset \Delta$ be an open topological disk and define $\rho_D$ to be the hyperbolic metric on $D$. Then $D$ is {\em round
  with respect to $f$} if $f(D)$ is round in $\chat$. We then define 
$$\rho_f(z) = \underset{D}{\inf} \rho_D(z).$$
where $D$ ranges over all round disks  containing $z$. By the Schwarz lemma, if $\rho_\Delta$ is the hyperbolic metric on $\Delta$ then $\rho_\Delta\le\rho_D$ for all disks $D$ contained in $\Delta$ with equality if and only if $D=\Delta$. Therefore $\rho_\Delta \le \rho_f$ with equality if and only if $f$ is the restriction of an element of $\psl$. In particular $\rho_f>0$. Upper bounds for $\rho_f$ are more subtle. The following theorem of Anderson will be a key tool for what follows.
\begin{theorem}{(Anderson, \cite[Thm. 4.2]{anderson-thesis})}\label{thurston_metric_bound}
$\rho_f(z) \le \rho_\Delta(z) \sqrt{1 + 2\|Sf\|_\infty}$.
\end{theorem}
Here $Sf$ is the {\em Schwarzian derivative} quadratic differential on $\Delta$ given by
$$Sf = \left(\left(\frac{f''}{f'}\right)' - \frac{1}{2} \left(\frac{f''}{f'}\right)^2\right)dz^2.$$
Then $\|Sf(z)\| = |Sf(z)|/\rho^2_\Delta(z)$ is a function on $\Delta$. In particular for any conformal automorphism $\gamma$ of $\Delta$ we have $\|S(f\circ \gamma)(z)\| = \|(Sf)(\gamma z)\|$. Furthermore the sup norm is given by
$$\|Sf\|_\infty = \sup_{z\in \Delta}\|Sf(z)\|.$$

\subsection{Epstein surfaces}
Using that $\chat$ can be naturally identified as the boundary of $\hthree$, we describe a construction of Epstein that associates a surface in $\hthree$ to a locally univalent map $f\colon\Delta\to\chat$ and a conformal metric $\rho$ on $\Delta$. 
	
Given a point $x \in \hthree$ let $\rho_x$ be the {\em visual metric} on $\chat$ centered at $x$. There are several ways to define $\rho_x$. We choose one that fits our needs for later. For $z \in \chat$ let $r$ be the geodesic ray based at $x$ that limits to $z$ at infinity. Then there will be a unique totally geodesic copy of $\htwo \subset \hthree$ that contains $x$ and is orthogonal to $r$. The hyperbolic plane will limit to a round circle in $\chat$. Let $D$ be the disk bounded by this circle that contains $z$ and $\rho_D$ its hyperbolic metric. We then define $\rho_x(z) = \rho_D(z)$. Note that $\rho_x$ is invariant under any isometry of $\hthree$ that fixes $x$. In fact, up to a normalization, this last property also determines $\rho_x$.

Given a conformal metric $\rho$ on a domain in $\chat$ containing a point $z$, we observe that the set $\ho_{\rho,z} = \{x \in \hthree | \rho_x(z) = \rho(z)\}$ is a horosphere. We will be interested in the horospheres associated to the push-forward metric $f_*\rho$. Unfortunately, as $f$ is only locally univalent, this push-forward is in general not well defined. To get around this, we define $f_*\rho(z)$ by restricting $f$ to a neighborhood of $z$ where $f$ is injective, pushing the metric forward on this neighborhood and then evaluating at $f(z)$. 

Let $T^1\hthree$ be the unit tangent bundle of $\hthree$ and $\pi\colon T^1\hthree\to\hthree$ the projection to $\hthree$. If $\rho$ is smooth Epstein shows that there is a unique smooth immersion $\widetilde\Ep_\rho\colon \Delta\to T^1\hthree$ such that $\widetilde\Ep_\rho(z)$ is an inward pointing normal to the horosphere $\ho_{f_*\rho, z}$ and when $\Ep_\rho = \pi\circ\widetilde\Ep_\rho$ is an immersion at $z$, the surface will be tangent to $\ho_{f_*\rho,z}$.  We emphasize that if $\rho$ is smooth then while $\widetilde\Ep_\rho$ will always be an immersion, $\Ep_\rho$ may not be. For example if $\rho_x$ is the visual metric for a point $x\in \hthree$ then $\widetilde\Ep_{\rho_x}$ is a diffeomorphism onto $T^1_x\hthree$ but $\Ep_{\rho_x}$ will be be the constant map to $x$.

The maps $\widetilde\Ep_\rho, \Ep_\rho$ have some nice properties. 
\begin{prop}{(Epstein \cite{epstein-envelopes})}\label{normal_flow}
Let $g_t\colon T^1\hthree\to T^1\hthree$ be the geodesic flow. Then $g_t\circ\widetilde\Ep_\rho = \widetilde\Ep_{e^t\rho}$. Furthermore if $\rho$ is smooth there are functions $\kappa^0_t, \kappa^1_t\colon \Delta\to (\R\backslash\{-1\})\cup\infty$ satisfying
$$\kappa^i_t(z) = \frac{\kappa_0^i(z) \cosh t + \sinh t}{\kappa_0^i(z) \sinh t + \cosh t}$$
such that if neither $\kappa^0_t(z), \kappa^1_t(z)$ are infinite then $\Ep_t$ is an immersion at $z$ and $\kappa^0_t(z), \kappa^1_t(z)$ are the principal curvatures. In particular if $t\ge \log\sqrt{\left| 1+\kappa_0^{i}(z)\right|/\left|1-\kappa_0^{i}(z)\right|}$ for $i=0,1$ then $\Ep_{e^t\rho}$ is an immersion and locally convex at $z$.
\end{prop}

The map $\Ep_\rho: \Delta \rightarrow \hthree$ is  the (parametrized) {\em Epstein surface} of $\rho$ associated to the locally univalent map $f$. We will be particularly interested in the Epstein surface $\Ep_{\rho_\Delta}$ associated to the hyperbolic metric $\rho_{\Delta}$ in $\Delta$. The importance of the Schwarzian derivative in studying the Epstein surface for the hyperbolic metric is evident in the following theorem.

\begin{theorem}{(Epstein \cite{epstein-envelopes})}\label{schwarzian_curvature}
The principal curvatures of the Epstein surface $\Ep_{\rho_\Delta}$ at the image of $z\in \Delta$ are $\frac{-\|Sf(z)\|}{\|Sf(z)\| \pm 1}$.
\end{theorem}

Theorem \ref{thurston_metric_bound} will follow from the following proposition.
\begin{prop}\label{convex_bound}
If $\rho$ is a smooth conformal metric and $\Ep_\rho$ is locally convex then $\rho_f \le \rho$.
\end{prop}

{\bf Proof:} Define a map $F\colon\Delta\times[0,\infty] \rightarrow \hthree$ by $F(z,t) = \Ep_{e^t\rho}(z)$ if $t\in[0,\infty)$ and $F(z,\infty) = f(z)$. By Proposition \ref{normal_flow}, $F$ restricted to $\Delta\times[0,\infty)$ will be an immersion to $\hthree$ and will extend continuously on $\Delta\times[0,\infty]$ to a map to $\hthree\cup\chat$. Since $F$ is an immersion, $F$ pulls back a hyperbolic structure on $\Delta\times[0,\infty)$ that is foliated by the Epstein surfaces. By convexity, a hyperbolic plane tangent to any Epstein surface in $\Delta\times [0,\infty)$ will be embedded and extend to a round disk on $\Delta = \Delta \times \{\infty\}$ with respect to $f$. For a point $z \in \Delta$, let $D$ be the round disk bounded by the boundary of the hyperbolic plane tangent to the Epstein surface at $(z,0)$. By definition, $\rho_f \le \rho_D$. On the other hand $\rho_D = \rho$ from the definition of the Epstein surface (and our normalization of the visual metric) and therefore $\rho_f(z) \le \rho(z)$ for all $z \in \Delta$.
\eproof

{\bf Proof of Theorem \ref{thurston_metric_bound}:}  By Theorem \ref{schwarzian_curvature}, the principal curvatures of $\Ep_{\rho_\Delta}$ at $\Ep_{\rho_\Delta}(z)$ are $\frac{-\|Sf(z)\|}{\|Sf(z)\|\pm 1}$. By the curvature equations in Proposition \ref{normal_flow} if $t> \sqrt{1+2\|Sf(z)\|}$ then the principal curvatures of $\Ep_{e^t\rho_\Delta}$ at $\Ep_{e^t\rho_\Delta}(z)$ are positive. So if $t>\sqrt{1+2\|Sf\|_\infty}$, then $\Ep_{e^t\rho_\Delta}$ is locally convex. The theorem then follows from Proposition \ref{convex_bound}. \eproof

If $Sf$ has small norm on a large neighborhood of $z \in \Delta$ then we can get stronger bounds on the Thurston metric.
\begin{cor}\label{local_thurston_bound}
If $\|Sf(z)\| \le K$ for all $z \in B(z_0,r)$ then
$$\rho_f(z_0) \le \rho_\Delta(z_0) \sqrt{1+2K}\coth(r/2).$$
\end{cor}

{\bf Proof:} 
Let $B = B(z_0,r)$. By the Schwarz Lemma $\frac{|Sf(z)|}{\rho_{B}(z)^2} \le \frac{|Sf(z)|}{\rho_\Delta(z)^2}$ and therefore by Theorem \ref{thurston_metric_bound}, $$\rho_{f|_B}(z_0) \le \rho_{B}(z_0)\sqrt{1+2K},$$ where $\rho_{f|_B}$ is the projective metric for $f$ restricted to $B$. By the definition of the Thurston metric $\rho_f(z_0) \le \rho_{f|_B}(z_0)$ and an explicit calculation shows that $\rho_{B}(z_0) = \rho_\Delta(z_0)\coth(r/2)$. This gives the desired inequality. \eproof

\subsection{The Thurston parameterization}
Let $P(\Delta)$ be locally univalent maps $f\colon\Delta\to \chat$ with the equivalence  $f\sim g$ if $f =\phi\circ g$ for some $\phi \in \psl$. Thurston described a natural parameterization of $P(\Delta)$ by $\cM\cL(\htwo)$ the space of measure geodesic laminations on $\htwo$. We briefly review this construction.

A round disk $D \subset \chat$ shares a boundary with a hyperbolic plane $\htwo_D \subset \hthree$. Let $r_D\colon D\to\hthree$ be the nearest point projection to $\htwo_D$ and $\tilde{r}_D\colon D\to T^1\hthree$ be the normal vector to $\htwo_D$ at $r_D(z)$ pointing towards $D$. We can 
use these maps to define a version of the Epstein map for $\rho_f$. In particular define $\widetilde\Ep_{\rho_f}\colon\Delta\to T^1\hthree$ by $\widetilde\Ep_{\rho_f}(z) = \tilde r_{f(D)}(f(z))$ where $D$ is the unique round disk with respect to $f$ such that $\rho_D(z) = \rho_f(z)$ and let $\Ep_{\rho_f}(z) = \pi\circ\widetilde\Ep_{\rho_f}(z) = r_{f(D)}(f(z))$. (For the existence of this disk see \cite[Theorem 1.2.7]{KT-projective}.) We also define $\widetilde\Ep_{e^t\rho_f} = g_t\circ\widetilde\Ep_{\rho_f}$ and $\Ep_{e^t\rho_f} = \pi\circ\widetilde\Ep_{e^t\rho_f}$.

The image of $\Ep_{\rho_f}$ is a {\em locally convex pleated plane}. More precisely, let $\cM\cL(\htwo)$ be measured geodesic laminations on $\htwo$ and $\cM\cL_0(\htwo) \subset \cM\cL(\htwo)$ the subspace of laminations with finite support. That is $\lambda \in \cM\cL_0(\htwo)$ if it is the union of a finite collection of disjoint geodesics $\ell_i$ with positive weights $\theta_i$. Then $\lambda$ determines a continuous map $p_\lambda\colon\htwo\to\hthree$, unique up to post-composition with isometries of $\hthree$, that is an isometry on the complement of the support of $\lambda$ and is ``bent'' with angle $\theta_i$ at $\ell_i$. By continuity we can extend this construction to a general $\lambda \in\cM\cL(\htwo)$. An exposition of the following theorem of Thurston can be found in \cite{KT-projective}.

\begin{theorem}\label{thurston_param}
Given $f \in P(\Delta)$ there exists maps $c_f\colon\Delta\to\htwo$ and $p_f\colon\htwo\to\hthree$ and a lamination $\lambda_f$ such that $p_f$ is a locally, convex pleated surface pleated along $\lambda_f$, $\Ep_{\rho_f} = p_f \circ c_f$ and the map $f\mapsto \lambda_f$ is a homeomorphism from $P(\Delta)\to\cM\cL(\htwo)$. Furthermore the maps $c_f\colon(\Delta, \rho_f)\to \htwo$ and $\Ep_{\rho_f}\colon(\Delta,\rho_f)\to \hthree$ are 1-Lipschitz.
\end{theorem}

Combined with Theorem \ref{thurston_metric_bound} 
we have an immediate corollary:
\begin{cor}\label{projective_lipschitz}
Given $f \in P(\Delta)$ the Epstein map $\Ep_{\rho_f}\colon (\Delta, \rho_\Delta)\to \hthree$ is $\sqrt{1+2\|Sf\|_\infty}$-Lipschitz.
\end{cor}

\subsection{Projective structures}
A projective structure $\Sigma$ on a surface $S$ is an atlas of charts to $\chat$ with transition maps the restriction of M\"obius transformations. We let $P(S)$ be the space of projective structures on $S$. One way to construct a projective structure on $S$ is to take an $f\in P(\Delta)$ such that there exists a Fuchsian group $\Gamma$ with $S = \Delta/\Gamma$ and a representation $\sigma\colon\Gamma\to\psl$ with $f\circ \gamma = \sigma(\gamma)\circ f$ for all $\gamma \in \Gamma$.  In fact every projective structure on $S$ arises in this way.

This description of projective structures allow us to associate a number of objects to a given projective structure. First we observe that a projective structure determines a complex structure  $X$ on $S$ and we let $P(X) \subset P(S)$ be projective structures on $S$ with underlying complex structure $X$. Given $\Sigma \in P(X)$, the Schwarzian derivative $Sf$ of $f$ will descend to a holomorphic quadratic differential $\phi_\Sigma$ on $X$. The lamination $\lambda_f$ will be $\Gamma$-equivariant and descend to measured lamination $\lambda_\Sigma$ on $S$. The hyperbolic metric $\rho_\Delta$ and the projective metric $\rho_f$ on $\Delta$ will also descend to conformal metrics $\rho_X$ and $\rho_\Sigma$ on $X$. 

In the equivariant setting Corollary \ref{projective_lipschitz} becomes:
\begin{theorem}\label{thurston_metric_projective}
Given a projective structure $\Sigma \in P(X)$ we have
$$\rho_\Sigma(z) \le \rho_X(z)\sqrt{1+2\|\phi_\Sigma\|_\infty}.$$
If $\|\phi_\Sigma(z)\| \le K$ for all $z \in B(z_0,r)$ then
$$\rho_\Sigma(z_0) \le \rho_X(z_0) \sqrt{1+2K}\coth(r/2).$$
\label{localS}
\end{theorem}

If the measured lamination $\lambda_\Sigma$ has support a finite collection of closed geodesics $\gamma_1, \dots, \gamma_n$ with weights $\theta_1,\dots, \theta_n$ then the length of $\lambda_\Sigma$ is $L(\lambda_\Sigma) = \sum \theta_i \ell(\gamma_i)$ where $\ell(\gamma_i)$ is the hyperbolic length of $\gamma_i$. This length extends continuously to general measure laminations on $S$.

We have the following useful relationship between the area of the projective metric and the length of the bending lamination.
\begin{lemma}\label{area_bending}
Given a projective structure $\Sigma \in P(S)$ with bending lamination $\lambda_\Sigma \in \cM\cL(S)$ we have $\area(\rho_\Sigma) = L(\lambda_\Sigma)+2\pi|\chi(S)|$.
\end{lemma}

{\bf Proof:} Both the area of the projective metric and the length of the bending lamination vary continuously in $P(S)$. The set of projective structures whose bending laminations is supported on finitely many geodesics is dense in $P(S)$ and the formula $\area(\rho_\Sigma) = L(\lambda_\Sigma) +2\pi|\chi(S)|$ holds on such laminations by direct computation. The lemma follows.\eproof

This immediately leads to bounds on the length:
\begin{theorem}\label{infinity_norm_bending_bound}
If $\lambda_\Sigma$ is the bending lamination for a projective structure $\Sigma$ then 
$$L(\lambda_\Sigma)\le 4\pi|\chi(\Sigma)|\|\phi_\Sigma\|_\infty.$$
\end{theorem}

{\bf Proof:} Squaring the inequality from  Theorem  \ref{thurston_metric_projective} we get a bound on the area of the projective metric in terms of the area of the hyperbolic metric:
$$\area(\rho_\Sigma) \le (1+2\|\phi_\Sigma\|_\infty) \area(\rho_X).$$
Subtracting $\area(\rho_X) = 2\pi|\chi(\Sigma)|$ from both sides and applying Lemma \ref{area_bending} we have
$$L(\lambda_\Sigma)\le 4\pi|\chi(\Sigma)|\|\phi_\Sigma\|_\infty$$
as claimed.\eproof

\subsection{Schwarzian bounds}
We recall the classical Nehari bound on the Schwarzian derivative. (The upper bound was proved independently by Kraus.)
\begin{theorem}{(Nehari \cite{nehari})}\label{kn_bound}
\begin{itemize}
\item If $f\colon\Delta\to\chat$ is univalent then $\|Sf\|_\infty \le \frac32$.
\item If $\|Sf\|_\infty \le \frac12$ then $f$ is univalent.
\end{itemize}
\end{theorem}
In particular, if $\Omega \subset \chat$ is a simply connected hyperbolic domain then the above theorem bounds the Schwarzian derivative of the uniformizing map $f\colon\Delta\to\Omega$. If $\Omega$ is hyperbolic but not necessarily simply connected we can still bound the Schwarzian for $f$ (which in this case will be a covering map) but our bounds depend on the injectivity radius of the hyperbolic metric of $\Omega$. Let $\inj_\Omega(z)$ be the supremum of the radii of embedded disks in $\Omega$ centered at $z$ and let
$$\delta_\Omega = \underset{z\in\Omega}{\inf} \inj_\Omega(z).$$
The following result bounding the Schwarzian in terms of $\inj_\Omega$ and $\delta_\Omega$ is due to Kra-Maskit. We include the short proof:
\begin{cor}{(Kra-Maksit, \cite[Lemma 5.1]{KM})}\label{compressible_bound}
Let $\Omega$ be a hyperbolic domain in the plane that is not simply connected and $f\colon\Delta\to\Omega$ be the uniformizing covering map. Then $\|Sf(z)\| \le \frac32\coth^2(\inj_\Omega(z)/2)$ and $\frac 12 \coth^2(\delta_\Omega/2) \le \|Sf\|_\infty$.
\end{cor}

{\bf Proof:} For each $z \in \Delta$ the restriction of $f$ to the disk $B=B(z,\inj_\Omega(z))$ is univalent. Applying Theorem \ref{kn_bound} we have that $$\frac{|Sf(z)|}{\rho^2_B(z)} \le \frac32$$ where $\rho_B$ is the hyperbolic metric on $B$. We also have $\rho_\Delta(z) = \tanh(\inj_\Omega(z)/2)\rho_B(z)$. The upper bound follows. 

Given any $\delta'> \delta_\Omega$ there exists a disk $B=B(z, \delta)$ such that $f|_B$ is not injective. Therefore, by Theorem \ref{kn_bound}, there exists a $z'\in B$ such that $$\frac{|Sf(z')|}{\rho^2_B(z')} \ge \frac12.$$ A calculation shows that $$\frac{\rho_B(z')}{\rho_\Delta(z')} \ge\frac{\rho_B(z)}{\rho_\Delta(z)} = \coth(\delta'/2)$$ so $\|Sf(z')\| \ge \frac12 \coth^2(\delta'/2)$. As this holds for all $\delta'>\delta_\Omega$ the lower bounds follows. \eproof

We  will only use the upper bound in what follows.


\subsection{Lipschitz maps and hyperbolic domes}
Let $\Omega \subset \chat$ be a hyperbolic domain and let $\Lambda = \chat\backslash\Omega$. The the {\em convex hull}, $H(\Lambda) \subset\hthree$ is the smallest closed convex subset of $\hthree$ whose closure in $\chat$ is $\Lambda$. The boundary of $H(\Lambda)$ is the {\em dome} of $\Omega$ which we denote $\dome(\Omega)$. With its intrinsic path metric, $\dome(\Omega)$ is a hyperbolic surface. The nearest point retraction of $\hthree$ to $H(\Lambda)$ extends to a continuous map $r\colon\Omega\to\dome(\Omega)$. We are interested in comparing the hyperbolic metric on $\Omega$ with the intrinsic path metric on $\dome(\Omega)$. 

We would like to relate the retract $r$ to an Epstein map. Let $f\colon\Delta\to\Omega$ be the uniformizing map. Then $f$ is a covering map and for any conformal metric $\rho$ on $\Delta$ that is invariant with respect to the covering, the Epstein map for $\rho$ will descend to a map with domain $\Omega$ which (in abuse of notation) we will continue to denote $\Ep_\rho\colon\Omega\to \hthree$. We then have
\begin{prop}\label{dome_epstein}
If $f\colon\Delta\to\Omega$ is the uniformizing map then $r\circ f = \Ep_{\rho_f}$.
\end{prop}

{\bf Proof:} Given $z \in \Omega$ there is a unique horosphere $\ho$ based at $z$ that intersects $H(\Lambda)$ at exactly one point with this point being the projection $r(z)$.

The hyperbolic plane tangent to $\ho$ at $r(z)$ is a support plane for $H(\Lambda)$ and it boundary bounds a round disk $D_z \subset \Omega$ which contains $z$. If $\rho_{D_z}(z) = \rho_f(z)$ then $\Ep_{\rho_f}(z) = r\circ f$ from the construction of the Epstein map for the projective metric.

By the definition of the projective metric $\rho_{D_z}(z) \ge \rho_f(z)$ so we just need to show that $\rho_{D_z}(z) \le \rho_f(z)$.

If $\rho_{D_z}(z) > \rho_f(z)$ then there exists a round disk $D \subset \Omega$ with $\rho_{D_z}(z) > \rho_D(z)$. Let $\ho'$ be the horosphere of points whose visual metrics agree with $\rho_D$ at $z$. Since $\rho_{D_z}(z) > \rho_D(z)$ the horosphere $\ho'$ bounds a horoball whose interior contains $\ho$. The open hyperbolic half space bounded by $D$ will contain the interior of this horoball and hence $\ho$. Since $\ho$ intersects $H(\Lambda)$, this open half space will intersect $H(\Lambda)$, a contradiction. \eproof

Combining this proposition with Theorem \ref{kn_bound} and Corollaries \ref{projective_lipschitz} and \ref{compressible_bound} we have the following:
\begin{theorem}\label{dome_lipschitz}
If $f\colon\Delta\to\Omega\subset \chat$ is a conformal homeomorphism then the retract $r\colon\Omega\to \dome(\Omega)$ is a $\sqrt{1+2\|Sf\|_\infty}$-Lipschitz map from the hyperbolic metric on $\Omega$ to the path metric on $\dome(\Omega)$. In particular, if $\Omega$ is simply connected then $r$ is 2-Lipschitz and if $\Omega$ is not simply connected with $\delta_\Omega >0$ then $r$ is $\sqrt{1+3\coth^2(\delta_\Omega/2)}$-Lipschitz.
\end{theorem}
When $\Omega$ is simply connected, Epstein-Marden-Markovic proved that the retract map was 2-Lipschitz \cite[Theorem 3.1]{EMM1}. When $\Omega$ is not simply connected Bridgeman-Canary showed that $r$ was $\left(A + \frac B{\delta_\Omega}\right)\mbox{-Lipschitz}$ for universal constants $A,B>0$ \cite[Corollary 1.8]{BC10}. Our bounds are better both when $\delta_\Omega$ is small and large. The simplicity of the proof here indicates one strength of our methods.

\subsection{Hyperbolic 3-manifolds}
The above result in turn can also be interpreted in terms of hyperbolic 3-manifolds. Let $\Gamma$ be a discrete, torsion free subgroup of $\psl$.  Let $\Omega$ be a component the domain of discontinuity of $\Gamma$ and let $\Gamma_\Omega \subset \Gamma$ be the subgroup that stabilizes $\Omega$. Then the projective structure $\Sigma = \Omega/\Gamma_\Omega$ is a component of the conformal boundary of the hyperbolic 3-manifold $M = \hthree/\Gamma$ and $X=\dome(\Omega)/\Gamma_\Omega$ is the component of the boundary of the convex core, $C(M)$, of $M$ that faces $\Sigma$. The nearest point retraction $M\to C(M)$ extends continuously to a map $r\colon \Sigma\to X$. Note that $\Sigma$ is incompressible in $M$ if and only if $\Omega$ is simply connected. If $\Sigma$ is compressible in $M$ then compressible curves in $\Sigma$ lift to homotopically non-trivial closed curves in $\Omega$. In particular the length of the shortest compressible geodesic in $\Sigma$ will be twice the injectivity radius of $\Omega$. In this setting Theorem \ref{dome_lipschitz} becomes:
\begin{theorem}\label{convex_core_lipschitz}
Let $M$ be a complete hyperbolic 3-manifold, $\Sigma$ the complex projective structure on a component of the conformal boundary of $M$ and $X$ the component of the boundary of the convex core of $M$ facing $\Sigma$. Then the retraction $r\colon\Sigma\to X$ is a $\sqrt{1+2\|Sf\|_\infty}$-Lipschitz map  from the hyperbolic metric on $\Sigma$ to the path metric on $X$. In particular, if $X$ is incompressible in $M$ then $r$ is 2-Lipschitz and if the length of every compressible curve on $X$ has length $\ge \delta>0$ then $r$ is $\sqrt{1+3\coth^2(\delta/4)}$-Lipschitz.
\end{theorem}

We can also apply the Schwarzian bounds to obtain bounds on the length of the bending lamination. In particular, Theorem \ref{infinity_norm_bending_bound} becomes:
\begin{theorem}\label{inj_bending_bound}
Let $\Sigma$ be a component of the projective boundary of a hyperbolic 3-manifold $M$ with bending lamination $\lambda_\Sigma$.  Then $L(\lambda_\Sigma)\le 4\pi |\chi(\Sigma)|\|\phi_\Sigma\|_\infty $. In particular
\begin{itemize}
\item If $\Sigma$ is incompressible then $L(\lambda_\Sigma)\leq 6\pi|\chi(\Sigma)|$.
\item If $\Sigma$ is compressible and the length of the shortest compressible curve is $\delta>0$ then $L(\lambda_\Sigma) \le 6\pi|\chi(\Sigma)|\coth^2(\delta/4)$.
\end{itemize}
\end{theorem}
The bound in the incompressible case was first obtained by Bridgeman-Canary in  \cite{BCrenorm}.  In the compressible case the bound in \cite{BC05} is $\left(\frac A\delta +B\right)|\chi(\Sigma)|$ which is stronger than the bound here. With more work our methods can obtain similar bounds as in \cite{BC05}. The proof is technical and this result won't be used in the rest of the paper.

\begin{theorem}\label{strong_compressible_bound}
If $\Sigma$ is a compressible component of the boundary of a hyperbolic 3-manifold with bending lamination $\lambda_\Sigma$ and the length shortest compressible curve is $\delta>0$ then $L(\lambda_\Sigma) \le \left(\frac{A}\delta +B\right)|\chi(\Sigma)|$ for universal constants $A,B>0$.
\end{theorem}

{\bf Proof:} The central idea is that that the ratio between the projective metric and the hyperbolic metric can only be large in the (compressible) thin part of the surface.

The complex projective structure $\Sigma$ is the quotient of a domain $\Omega \subset \chat$.
Let $X$ be the conformal structure on $\Sigma$ with hyperbolic metric $\rho_X$. Then $\Omega$ is a covering space of $X$ and the hyperbolic metric $\rho_\Omega$ is the lift of $\rho_X$. Similarly the Schwarzian $\phi_\Omega$ on $\Omega$ is the lift of the Schwarzian $\phi_\Sigma$ on $\Sigma$ and $\rho_{\tilde\Sigma}$ is the lift of the projective metric $\rho_\Sigma$. We would like to bound above the ratio $\rho_\Sigma(z)/\rho_X(z)$. To do this we will use Corollary \ref{projective_lipschitz} which will require us to bound the Schwarzian $\|\phi_\Sigma(z')\|$ for all $z'$ in the disk $B(z,1)$. (The choice of radius $1$ is essentially arbitrary.) To bound $\phi_\Sigma$ we will  use Corollary \ref{compressible_bound} to bound $\phi_\Omega$ and then use that $\phi_\Omega$ is the lift of $\phi_\Sigma$.

For $z\in\Omega$ let $$\inj^1_\Omega(z) = \underset{z'\in B(z,1)}{\inf} \inj_\Omega(z').$$ A simple estimate gives that there exists a constant $A_0>0$ such that $\inj_\Omega^1(z) \ge \inj_\Omega(z)/A_0$. (This holds for any complete hyperbolic surface.) 

Let $\epsilon>0$ be the two-dimensional Margulis constant and let $C \subset \Omega$ be a component of $\epsilon$-thin part $\Omega^{<\epsilon} = \{z \in \Omega\  |\  \inj_\Omega(z) < \epsilon\}.$ There is also a constant $A_1>0$ such that $\frac 32\coth^2(x/2) \le \frac {A_1}{x^2}$ for $x \leq \epsilon$.
Then by Corollary \ref{compressible_bound}, for all  $z'\in B(z,1)$ with $z \in C$ we have $$\|\phi_\Omega(z')\| \le \frac {A_1}{\inj^1_\Omega(z)^2} \le \frac{A_0^2A_1}{\inj_\Omega(z)^2}.$$
Applying Corollary \ref{local_thurston_bound} we see that for $z \in C$
$$\rho^2_{\tilde\Sigma}(z) \le \rho^2_\Omega(z)\left(1+\frac{2A_0^2A_1}{\inj_\Omega(z)^2}\right)\coth^2(1/2).$$

We want to bound the area of $C$ in the projective metric. We let $\ell$ be the length of the core geodesic of $C$ in the hyperbolic metric. We give $C$ coordinates $S^1\times (-w(\ell), w(\ell))$ where $S^1 \times \{0\}$ is the geodesic and is parameterized by arc length and each $\{\theta\}\times (-w(\ell), w(\ell))$ is a geodesic segment orthogonal to the core geodesic. The area form for the hyperbolic metric is then $\cosh t d\theta dt$. The constant $w(\ell)$ is chosen such that $\inj_\Omega(\theta, \pm w(\ell)) = \epsilon$. Another basic estimate in hyperbolic geometry gives that there exists $A_2>1$ such that
$$\ell e^{|t|}/A_2 \le \inj_\Omega(\theta, t).$$
Here it is important that $z =(\theta,t) \in C$ is in the $\epsilon$-thin part.

We now calculate the area of $C$ in the projective metric:
\begin{eqnarray*}
\area(\rho_{\tilde\Sigma}|_C) & = & \int_{-w}^w\int_0^\ell  \cosh t\rho^2_{\tilde\Sigma}/\rho^2_\Omega  d\theta dt \\
& \le & \coth^2(1/2)\left(\area(\rho_{\Omega|_C}) + \int_{-w}^w\int_0^\ell  \frac{2A_0^2A_1}{\inj_\Omega(\theta,t)^2}\cosh td\theta dt\right).
\end{eqnarray*}
We use the lower bound on the injectivity radius to bound the remaining integral:
\begin{eqnarray*}
\int_{-w}^w\int_0^\ell  \frac{2A_0^2A_1}{\inj_\Omega(\theta,t)^2}\cosh td\theta dt & \le & 2\int_{0}^w\int_0^\ell  \frac{2A_0^2A_1A_2^2}{\ell^2 e^{2t}}\cosh td\theta dt\\
& = & 2\int_{0}^w \frac{2\ell A_0^2A_1A^2_2}{\ell^2e^{2t}}\cosh td\theta dt\\
& \le & A_3 \int_0^w \frac{e^t}{\ell e^{2t}} dt \\
& = &\frac{A_3}\ell  (1 - e^{-w}) \le \frac{2A_3}\ell .
\end{eqnarray*}
Since $\ell  \ge \delta$ this becomes
$$\area(\rho_{\tilde\Sigma}|_C) \le \coth^2(1/2) \left(\area(\rho_\Omega|_C) + \frac{A_3}{\delta}\right).$$

Given a point $z \in \Sigma$ let $\tilde z \in \Omega$ be a point in the pre-image of $z$. We then define $\tilde\inj_X(z) = \inj_\Omega(\tilde z)$ and observe that this definition is independent of our choice of $\tilde z$. Injectivity radius can only increase in a cover so $\inj_X(z) \le \tilde\inj_X(z)$. The {\em compressible $\epsilon$-thin part} is the set of points
$$X^{<\epsilon}_c=\{z \in X\  |\  \tilde\inj_X(z) < \epsilon\}.$$
If $C$ is a component of the compressible $\epsilon$-thin part then each component of the pre-image of $C$ in $\Omega$ will be contained in a component $\tilde C$ the $\epsilon$-thin part of $\Omega$ and we'll have $\area(\rho_{\Sigma}|_C) \le \area(\rho_{\tilde\Sigma}|_C)$. Furthermore each $C$ will contain a simple closed geodesic so there can be at most $3g-3 = \frac 32|\chi(\Sigma)|$ components of the compressible $\epsilon$-thin part and therefore
$$\area(\rho_\Sigma|_{X^{<\epsilon}_c}) \le \coth^2(1/2)\left(\area(\rho_X|_{X^{<\epsilon}_c}) + \frac{3A_3}{2\delta}|\chi(\Sigma)|\right).$$

On the other hand if $\tilde\inj_X(z) \ge \epsilon$ then for $z' \in B(z,1)$ we have as above that $\tilde\inj_X(z') \ge \epsilon/A_0$. Therefore by Corollary \ref{local_thurston_bound} 
$$\rho^2_\Sigma(z)/\rho^2_X(z) = \rho^2_{\tilde\Sigma}(\tilde z)/\rho^2_\Omega(\tilde z) \le \left(1 +3\coth^2\left(\frac{\epsilon}{2A_0}\right)\right) \coth^2(1/2)= A_4.$$
Therefore we have that 
$$\area(\rho_{\Sigma}|_{X^{\ge \epsilon}_c}) \le A_4 \area(\rho_X|_{X^{\ge \epsilon}_c})$$
where $X^{\ge \epsilon}_c$ is the compressible $\epsilon$-thick part of $X$.

Letting $A= \frac32A_3\coth^2(1/2)$ and $B=2\pi A_4$ and combining our two area bounds we have
$$L(\lambda_\Sigma) \le \area(\rho_\Sigma) = \area(\rho_\Sigma|_{X^{<\epsilon}_x}) + \area(\rho_\Sigma|_{X^{\ge \epsilon}_c}) \le |\chi(\Sigma)|\left(\frac A{\delta} + B\right).$$\eproof
\subsection{$L^2$-bounds for the bending lamination}

Given a quadratic differential $\phi$ on hyperbolic surface $X$ with metric $\rho_X$, the ratio $|\phi|/\rho_X^2$ is a function on $X$. We define the $L^2$-norm of $\phi$ to be the $L^2$-norm of this function with respect to the hyperbolic metric.
In order to prove our main theorem about renormalized volume, we will need a bound on $L(\lambda_\Sigma)$ in terms of the $L^2$-norm of the quadratic differential $\phi_\Sigma$. We begin with the following lemma.
  
\begin{lemma}
\label{L2injectivity_bound}
Let $\phi$ be a holomorphic quadratic differential on a hyperbolic surface $X$. 

Then
$$\|\phi\|_{2} \ge 2\sqrt\frac\pi3 \tanh^2(\inj_X(z)/2)\|\phi(z)\|.$$
\end{lemma}

{\bf Proof:} Let $B=B(z,r)$ be the disk centered at $z$ of radius $r=\inj_X(z)$. Let $\|\phi\|_{X,2}$ be the $L^2$-norm of $\phi$ on $X$ and  $\|\phi\|_{B,2}$ be the $L^2$-norm of $\phi$ on $B$. Then $\|\phi\|_{X,2} \ge \|\phi\|_{B,2}$ by the Schwarz lemma. By \cite[Lemma 5.1]{bromberg}  we have $\|\phi\|_{B,2} \ge 2\sqrt{\pi/3}\|\phi(z)\|_B$ where $\|\phi(z)\|_B$ is the norm of $\phi$ on $B$. Comparing the complete hyperbolic metric on $B$ to that on $X$ we see $\|\phi(z)\|_B = \tanh^2(r/2)\|\phi(z)\|_X$. \eproof

We now combine the above with the prior results to obtain comparisons of the Thurston metric and Poincar\'e metric for quadratic differentials with small $L^2$-norm on the thick part of the surface. For $\epsilon > 0$ we define the $\epsilon$ thick-thin decomposition to be $X^{\geq \epsilon} = \{ z \in X\ |\ \inj_X(z) \geq \epsilon\}$ and $X^{< \epsilon} = \{ z \in X\ |\ \inj_X(z) < \epsilon\}$.

\begin{lemma}
Let $\Sigma \in P(X)$  be a projective structure  such that $||\phi_\Sigma||_2 \leq \epsilon^5$. Then for $z \in X^{\geq \epsilon}$
$$\rho_\Sigma(z) \leq (1+F(\epsilon))\rho_X(z)$$
where $F(t) \simeq (2 + 4\sqrt{3/\pi})t$ as $t \rightarrow 0$. 
\end{lemma}

{\bf Proof:} We can assume that $\epsilon<1$ and then define $r > 0$ such that $\epsilon =e^{-r}$.
Let $z \in X^{\geq \epsilon}$. For $w \in B(z,r)$ then a simple calculation shows that $\inj_X(w) \geq \inj_X(z) e^{-r} \geq \epsilon^2$.
 This follows from the fact that for $C$ a hyperbolic annulus with core geodesic of length $\ell$ then 
$$\sinh(\inj_C(x)) = \sinh(\ell/2)\cosh(d(x))$$
where $d(x)$ is the distance from $x$ to the geodesic. Comparing two points $x,y$ with $\inj_C(x) > \inj_C(y)$ one obtains 
$$\frac{\inj_C(x)}{\inj_C(y)} \leq  \frac{\sinh(\inj_C(x))}{\sinh(\inj_C(y))} \leq \frac{\cosh(d(x))}{\cosh(d(y))}  \leq e^{d(x)-d(y)}.$$

 Therefore for  $w \in B(z,r)$ by Lemma \ref{L2injectivity_bound}
$$\|\phi_\Sigma(w)\| \leq  \sqrt{\frac{3}{4\pi}}\left(\frac{\|\phi_\Sigma\|_{2}}{\tanh^2(\epsilon^2/2)}\right) \leq \sqrt{\frac{3}{4\pi}}\left(\frac{\epsilon^5}{\tanh^2(\epsilon^2/2)}\right)$$
 Therefore by the local bound in Theorem \ref{localS} we have
\begin{eqnarray*}
\frac{\rho_\Sigma(z)}{\rho_X(z)} &\leq &\sqrt{1+\sqrt{\frac{3}{\pi}}\left(\frac{\epsilon^5}{\tanh^2(\epsilon^2/2)}\right)}\coth(r/2) \\ &= & \sqrt{1+\sqrt{\frac{3}{\pi}}\left(\frac{\epsilon^5}{\tanh^2(\epsilon^2/2)}\right)}\left(\frac{1+\epsilon}{1-\epsilon}\right) = 1+F(\epsilon).
\end{eqnarray*}
Computing the first two terms of the Taylor series shows that as $t \rightarrow 0$
 $$F(t) \simeq \left(2+4\sqrt{\frac{3}{\pi}}\right) t.$$ 
 \eproof

We now use the above to get prove the $L^2$-bound on the length of the bending lamination.

\smallskip
\begin{theorem}\label{L2_bend_bound}
Let $\Sigma \in P(X)$  be a projective structure with Schwarzian quadratic differential $\phi_\Sigma$ with $\|\phi_\Sigma\|_\infty \le K$.  Then 
$$L(\lambda_{\Sigma}) \leq  2\pi |\chi(X)| G_K(||\phi_{\Sigma}||_2) $$
where $G_K(t)\sim t^{1/5}$ as $t \rightarrow 0$.
\end{theorem}

{\bf Proof:}
We let $\epsilon = \|\phi\|^{1/5}_2$. As $\|\phi_\Sigma\|_\infty \le K$, by Theorem \ref{thurston_metric_projective} we have $\rho_\Sigma(z) \leq \sqrt{1+2K}\rho_X(z)$ for all $z$. We decompose $X$ into the thick-thin pieces.
$$\area(\rho_\Sigma) = \int_{X^{\geq\epsilon}} \rho_\Sigma^2 + \int_{X^{<\epsilon}}  \rho_\Sigma^2 \leq \int_{X^{\geq\epsilon}} (1+F(\epsilon))^2\rho_X^2  +  \int_{X^{<\epsilon}} (1+2K)\rho_X^2.$$
Therefore
$$\area(\rho_\Sigma) \leq (1+F(\epsilon))^2 \area(\rho_{X^{\geq\epsilon}}) + (1+2K)\area(\rho_{X^{<\epsilon}}).$$
Since $\area(\rho_{X^{\geq\epsilon}}) \leq \area(\rho_X) = 2\pi|\chi(S)|$ and for the genus $g$ surface $S$ there are at most $(3g-3)$  $\epsilon$-thin parts each with area bounded by $2\epsilon$ we have 
$$2\pi|\chi(S)| + L(\lambda_{\Sigma}) \leq (1+F(\epsilon))^2 2\pi|\chi(S)| + (1+2K)(3g-3)2\epsilon.$$
Since $|\chi(S)| = 2g-2$ when we apply Lemma \ref{area_bending} we have
\begin{eqnarray*}
L(\lambda_{\Sigma}) &\leq & 2\pi|\chi(S)|\left((1+F(\epsilon))^2+\frac{3\epsilon}{2\pi}(1+K)-1\right)\\ & = &2\pi|\chi(S)| \left(2F(\epsilon)+F(\epsilon)^2+\frac{3\epsilon}{2\pi}(1+2K)\right).
\end{eqnarray*}
\eproof

\section{Renormalized Volume}
We now describe the renormalized volume for a convex cocompact hyperbolic 3-manifold $M$. We also review many of its fundamental properties as developed by Krasnov and Schlenker. While it will take some setup before we state the definition, we will see that  renormalized volume  has many nice properties that make its definition natural.

\subsection{The $W$-volume}
Throughout this subsection and the next, we fix a convex cocompact hyperbolic 3-manifold $M$ and let $\partial_c M$ be its conformal boundary, $\Sigma$ its projective boundary, $\lambda_M$ the bending lamination of the convex core and $\phi_M$ the Schwarzian derivative of $\Sigma$. We also let $\rho_M$ be the hyperbolic metric on $\partial_c M$ and $\rho_\Sigma$ the projective metric determined by $\Sigma$.

Let $N \subset M$ be a smooth, compact convex submanifold of $M$ with $C^{1,1}$ boundary. Here, and in what follows, $N$ is convex if every geodesic segment with endpoints in $N$ is contained in $N$. Then the {\em $W$-volume} of $N$ is
$$W(N) = \vol(N) - \frac{1}{2}\int_{\partial N} H da$$
where $H$ is the mean curvature function on $\partial N$. That is, $H$ is the average of the principal curvatures or, equivalently, one half the trace of the shape operator. The $C^{1,1}$ condition (the normal vector field is defined everywhere and is lipschitz) implies that $H$ is defined almost everywhere and that the integral
$$2\int_{\partial N} H da$$
is the variation of area of $\partial N$ under the normal flow.

We let $N_t$ be the $t$-neighborhood of $N$ in $M$. Then there is a very simple formula for the $W$-volume of $N_t$ in terms of $N$.
\begin{prop}{(Krasnov-Schlenker, \cite{KS08})}\label{scaling}
Let $M$ be a a convex cocompact hyperbolic 3-manifold and $N$ a convex submanifold with $C^{1,1}$ boundary. Then
$$W(N_t) = W(N) - t\pi \chi(\partial_c M).$$
\end{prop}

As defined, the $W$-volume is a function on the space of convex submanifolds of $M$ with $C^{1,1}$ boundary. We would like to reinterpret it as a function on smooth, conformal metrics on $\partial_c M$. We need the following lemma:

\begin{lemma}\label{convex_gives_epstein}
Let $H$ be a closed convex submanifold of $\hthree$ and let $\Lambda = \overline H \cap \chat$ and $\Omega = \chat \backslash \Lambda$. Then there exists a conformal metric $\rho=\rho_H$ on $\Omega$ such that $\Ep_\rho$ is the nearest point retraction $r\colon\Omega\to \partial H$. If $\gamma \in \psl$ with $\gamma(H) = H$ then $\gamma^* \rho = \rho$.

In particular, if $N\subset M$ is a convex submanifold of a convex cocompact hyperbolic manifold $M$ then there exists a smooth metric $\rho=\rho_N$ on $\partial_c M$ such that $\Ep_\rho = r$ where $r\colon \partial_c M\to \partial N$ is the nearest point retraction.

Finally if $N_t$ is the $t$-neighborhood of $N$ then $\rho_{N_t} = e^t \rho_N$.
\end{lemma}

{\bf Proof.} For each $z \in \Omega$ there is a unique horosphere $\ho_z$ based at $z$ that intersects $H$ at exactly one point and $r(z)$ is the point of intersection.  We then define $\rho(z) = \nu_{r(z)}$ where $\nu_{r(z)}$ is the visual metric. Then $r$ satisfies all the properties of the Epstein map for $\rho$ and since the Epstein map is unique we have $r= \Ep_\rho$. The construction is clearly equivariant. Equivariance, implies the second paragraph and the last statement then follows from Proposition \ref{normal_flow}.
\eproof

We then have a nice formula the integral of the mean curvature in terms of the of the area of $\rho_N$ and $\partial N$.

\begin{lemma}\label{equivalent_W}
Let $N$ be a smooth convex submanifold of a convex cocompact hyperbolic 3-manifold $M$. Then
$$\int_{\partial N} H da = \frac12\area(\rho_N) - \area(\partial N) - \pi \chi(\partial M).$$
Furthermore if $\rho_N = \rho_M$  then
$$\int_{\partial N} H da = ||\phi_M||^2_2$$.
\end{lemma}
{\bf Proof:} Let $B\colon T(\partial N)\to T(\partial N)$ be the shape operator given by $B(v) = -\nabla_v n$ where $n$ is the normal vector field to $\partial N$ . In particular, the eigenvalues of $B$ are the principal curvatures of $\partial N$. Then
$$H = \frac12 \tr(B) = \frac14\left(\det(I+B) - \det(I-B)\right)$$
where $I\colon T(\partial N)\to T(\partial N)$ is the identity operator. An elementary calculation shows that   the pullback via the retraction $r\colon \partial_cM \to \partial N$ of the 2-form $\det(I+B)da$ is the area form for the metric $\rho_N$ on $\partial_c M$ (see \cite[5.3]{KS08}).  Therefore
$$\area(\rho_N) = \int_{\partial N} \det(I+B)da.$$
On other hand
$$\det(I+B) + \det(I-B) = 2 +K$$
where $K = \det(B) - 1$ is the Gaussian curvature of $\partial N$. Therefore
$$\int_{\partial_N} (\det(I+B) + \det(I-B) )da = 2\area(\partial N) + 2\pi\chi(\partial N).$$
Rearranging terms proves the first statement in the lemma.

For $\rho_N = \rho_M$ the hyperbolic metric, by Theorem \ref{schwarzian_curvature}  the principal curvatures at $r(z)$ are $\frac{-\||\phi_M(z)||}{||\phi_M(z)||\pm 1}$. Therefore if $da^*$ is the area form for $\rho_M$ then
$$\area(\partial N) = \int_{\partial N} da = \int_{\partial_c M} \frac{1}{\det(I+B)}da^* = \int_{\partial_cM} (1-||\phi_M(z)||^2)da^* = \area(\rho_M) -||\phi_M||_2^2. $$
Therefore as $\area(\rho_M) = 2\pi|\chi(\partial M)|$ the result follows.
\eproof

This gives us an alternate way of defining the $W$-volume by setting
$$W(N) = \vol(N) - \frac14\area(\rho_N) + \frac12\area(\partial N) +\frac12 \pi \chi(\partial N).$$

Note that the definition makes sense even if the boundary $N$ is not $C^{1,1}$. Also, regardless of the regularity of $N$, the $t$-neighborhood $N_t$ will always have $C^{1,1}$ boundary. In particular, the scaling property, Proposition \ref{scaling}, still holds for this alternative definition of the $W$-volume even when the boundary of $N$ is not $C^{1,1}$. One advantage of this definition is that we can use to see that the $W$-volume varies continuously.

\begin{prop}\label{Hausdorff_continuity}
The $W$-volume is continuous on the space of compact convex submanifolds of $M$ with the Hausdorff topology. The map from compact, convex submanifolds to metrics on $\partial_c M$ is continuous in the $L^\infty$-topology.
\end{prop}

{\bf Proof:} Fix a convex submanifold $N$ and let $V_i$ be convex submanifolds such that the distance between $N$ and $V_i$ in the Hausdorff metric is less then $1/i$. We can assume that $N \subset V_i$ for if not we can replace $V_i$ with its $1/i$-neighborhood. By Proposition \ref{scaling}, the $W$-volume of the $V_i$ will converge to $W(N)$ if and only if the $W$-volume of the $1/i$-neighborhoods also converge to $W(N)$.

To see that the $W$-volume converges we first observe that volume is continuous in the Hausdorff topology on the space of convex submanifolds. Next we note that the nearest point retraction of $\partial V_i$ to $\partial N$ is 1-Lipschitz so $\area(\partial V_i) \ge \area(\partial N)$. Since $V_i \subset N_{1/i}$ we similarly have that $\area(\partial N_{1/i}) \ge \area(\partial V_i)$. We also have $\area(\partial N_{1/i}) \to \area(\partial N)$ and therefore $\area(\partial V_i) \to \area(\partial N)$.

To compare the metrics $\rho_N$ and $\rho_{V_i}$, fix a point $z \in \partial_c M$ and let $\ho_z$ and $\ho^i_z$ be the horospheres based at $z$ that meet $N$ and $V_i$, respectively, in a single point. Then $\ho^i_z$ will be disjoint from $N$ but its $1/i$-neighborhood will intersect $N$. This implies that $1 \le \rho_{V_i}(z)/\rho_N(z) \le e^{1/i}$. It follows that the map from convex submanifolds with the Hausdorff topology to the space of conformal metrics with the $L^\infty$-topology is continuous. Therefore $\area(\rho_N)$ varies continuously in $N$ and this, along with the previous paragraph, implies that the $W$-volume varies continuously.
\eproof

Let $\mathcal M(\partial_c M)$ be continuous conformal metrics on $\partial_c M$ with the $L^\infty$-topology and let $\mathcal M_C(\partial_c M)$ be the subspace of metrics $\rho$ such that  there exists a convex submanifold $N$ with $\rho_N = e^t \rho$ for some $t\in \mathbb R$. We can then define the $W$-volume as a function on $\mathcal M_C(\partial_c M)$ by setting
$$W(\rho) = W(N) +t\pi\chi(\partial_c M).$$
Note that $\mathcal M_C(\partial_c M)$ will not be all continuous metrics. For example a metric that locally has the form $\frac{|dz|}{1+|z|}$ will not be in $\mathcal M_C(\partial_c M)$. However we have the following:

\begin{prop}{(Krasnov-Schenker \cite[Theorem 5.8]{KS08})}\label{eventually_convex}
Let $\rho$ be a smooth metric on $\partial_c M$. Then for $t$ sufficiently large there exists a convex submanifold $N\subset M$ such that $e^t \rho = \rho_N$. In particular $\rho \in \mathcal M_C(\partial_c M)$.
\end{prop}

We are now finally in position to define the renormalized volume. We let $\rho_M$ be the hyperbolic metric on $\partial_cM$ and define
$$\Rvol(M) = W(\rho_M).$$
We have, by the Lemma \ref{equivalent_W}, that if $N$ is the submanifold corresponding to $\rho_M$ then
 $$V_R(M) = \vol(N) - ||\phi_M||^2_2$$

\begin{theorem}{(Schlenker \cite[Proposition 3.11, Corollary 3.8]{schlenker-qfvolume})}
Let $M$ be convex cocompact. Then
\begin{itemize}
\item (Maximality) If $\rho\in\mathcal M_C(\partial_c M)$ with $\area(\rho) = \area(\rho_M)$ then $W(\rho) \le \Rvol(M)$ with equality if and only if $\rho=\rho_M$.
\item (Monotonicity) If $\rho_0, \rho_1 \in \mathcal M_C(\partial_c M)$ have non-positive curvature on $\partial_c M$ and $\rho_0 \leq \rho_1$ then
$$W(\rho_0) \leq W(\rho_1).$$
\end{itemize}
\label{Wprops}
\end{theorem}

\subsection{Bounds on Renormalized volume}

For quasifuchsian manifolds Schlenker used the $W$-volume of the convex core to get upper and lower bounds on the renormalized volume (see (Schlenker, \cite[Theorem 1.1]{schlenker-qfvolume}). This  generalized easily to convex cocompact 3-manifolds with incompressible boundary (see \cite[Theorem 1.1]{BCrenorm}).  

Schlenker's approach was to use the monotonicity property and maximality property of W-volume. If $N$ is the convex core of $M$ then by Proposition \ref{dome_epstein}, the metric at infinity $\rho_N = \rho_\Sigma$, the projective metric. As the projective metric is non-positively curved and is greater than the hyperbolic metric  the monotonicity property of $W$-volume (Theorem \ref{Wprops}) gives an upper bound on the renormalized volume. By rescaling the metric so that it has the same area as the hyperbolic metric, the maximality property (again Theorem \ref{Wprops}) gives a lower bound on renormalized volume.

%
%
%

In order to obtain $L^2$-bounds, we will use the same strategy as Schlenker. We first need the following Theorem. The upper bound is due to Schlenker (see \cite{schlenker-qfvolume}) and the lower bound  is a simple application of the monotonicity and maximality properties of renormalized volume.

\begin{theorem}
Let $M$ be a convex cocompact hyperbolic 3-manifold  then
$$V_C(M) -\frac12L(\lambda_M) \leq V_R(M) \leq V_C(M) -\frac14L(\lambda_M)$$
and $V_C(M) = V_R(M)$ if and only if the convex core of $M$ has totally geodesic boundary.
\end{theorem}

{\bf Proof:} As noted above, the metric at infinity $\rho_N$ for the convex core $N$ is the projective metric $\rho_\Sigma$. Using our formula for the $W$-volume in terms of area we compute that 
$$W(\rho_\Sigma) =W(N)=\Cvol(M) - \frac14\area(\rho_\Sigma) + \frac12\area(\partial N) + \frac12\pi\chi(\partial N).$$
By Lemma \ref{area_bending}, $\area(\rho_\Sigma) = L(\lambda_M) + 2\pi|\chi(S)|$. Since the boundary of the convex core is a hyperbolic surface we have $\area(\partial N) = 2\pi|\chi(S)|$. Therefore
$$W(\rho_\Sigma) = \Cvol(M) -\frac14L(\lambda_M).$$

Since $\rho_\Sigma \ge \rho_M$, by the monotonicity property we have 
$W(\rho_\Sigma) \ge \Rvol(M)$ so 
$$\Cvol(M) \ge \Cvol(M) - \frac14L(\lambda_M) \ge \Rvol(M)$$ with $V_C(M) = V_R(M)$ if and only if $L(\lambda_M) = 0$. This proves the upper bound.

For the lower bound let $$\hat\rho_\Sigma = \sqrt{\frac{\area(\rho_\Sigma)}{\area(\rho_M)}}\rho_\Sigma.$$
Then $\area(\hat\rho_\Sigma) = \area(\rho_M)$ so by the maximality property (Theorem \ref{Wprops}) $W(\hat\rho_\Sigma) \le \Rvol(M)$. Similarly by the scaling property (Theorem \ref{Wprops}) and the formula for $\area(\rho_\Sigma)$ we have
$$W(\hat\rho_\Sigma) = W(\rho_\Sigma) -\frac\pi2\log\left(1 +\frac{L(\lambda_M)}{\area(\rho_M)}\right)|\chi(\partial M)|.$$
As $\log(1+x) \leq x$  and $\area(\rho_M) = 2\pi|\chi(\partial M)|$ we have
$$V_R(M) \geq W(\hat\rho_\Sigma) \geq W(\rho_\Sigma)  -\frac14L(\lambda_M) = V_C(M) -\frac12L(\lambda_M).$$
Thus it follows that
$$V_C(M) -\frac12L(\lambda_M) \leq V_R(M) \leq V_C(M) -\frac14L(\lambda_M).$$
We therefore have $V_C(M) = V_R(M)$ if and only if $L(\lambda_M) = 0$. Thus $V_C(M) = V_R(M)$ if and only if $M$ has totally geodesic boundary.
\eproof

\medskip

\noindent
Combining the $L^2$-bound for length in Theorem \ref{L2_bend_bound} with the above theorem we obtain the following;\medskip

\noindent
{\bf Theorem \ref{L2_squeeze}} {\em
There is a function $G(t) \sim t^{1/5}$ such that  if $M$ is a convex cocompact hyperbolic 3-manifold with incompressible boundary then
$$\Cvol(M) - |\chi(\partial M)|G(||\phi_M||_2) \leq \Rvol(M) \leq \Cvol(M)$$
and $\Rvol(M) = \Cvol(M)$ if and only if the convex core of $M$ has totally geodesic boundary.
}
\medskip

{\bf Proof:}
As $M$ has incompressible boundary, then we have the Nehari bound $||\phi_M|| \leq \frac32 $. 
By Theorem \ref{L2_bend_bound} we have
$$L(\lambda_M) \le 2\pi|\chi(\partial_c M)| G_{\frac32}(\|\phi_M\|_2)$$
where $G_{\frac32}(t) \sim t^{1/5}$ and the result follows with $G = \pi .G_{\frac32}$.
\eproof

The results here should be compared to  earlier work of Bridgeman-Canary. For manifolds with incompressible boundary they prove a lower bound where the function $G(t)$ in Theorem \ref{L2_squeeze} is replaced by a universal constant (see \cite[Theorem 1.1]{BCrenorm}). For manifolds with compressible boundary they give upper and lower bounds on the $V_R(M)$ that depend on the length of the shortest compressible curve in the boundary (see \cite[Theorem 1.3]{BCrenorm}). In particular one can produce a sequence of Schottky manifolds (convex cocompact hyperbolic structures on a handlebody) of fixed genus whose convex core volume is bounded above but the length of the shortest compressible curve approaches zero. Then the Bridgeman-Canary bounds imply that the renormalized volume of this sequence limits to $-\infty$.


\subsection{The gradient flow of $\Rvol$}
Let $N$ be a compact, hyperbolizable 3-manifold with incompressible boundary and let $CC(N)$ be the space of convex cocompact hyperbolic structures on $N$. Then for each $M\in CC(N)$ the map $M \mapsto \partial_c M$ defines an isomorphism from $CC(N)$ to $\Teich(\partial N)$. The renormalized volume is  then a function on $CC(N)$ and, via the above identification, a function on $\Teich(\partial N)$.  By \cite[Corollary 8.6]{KS08}, $\Rvol$ is a smooth function and we let $V$ be the gradient flow of $-\Rvol$ with respect to the Weil-Petersson metric on $\Teich(\partial N)$.

Recall that a tangent vector to Teichm\"uller space is given by a Beltrami differential and a cotangent vector is a holomorphic quadratic differential.
\begin{theorem}{(Krasnov-Schlenker, \cite{KS08})}\label{Rvolvariation}
Given $M\in CC(N)$ and $\mu \in T_{\partial_c M}\Teich(\partial N)$ we have
$$d\Rvol(\mu) = \Re \int_{\partial_c M} \mu\phi_M.$$
\end{theorem}
Using this variational formula, we get an explicit description of the gradient flow.
\begin{prop}
The flow for $V$ is defined for all time and for each $M \in CC(N)$, $V(M) = -\frac{\bar\phi_M}{\rho^2_M}$.
\end{prop}

{\bf Proof:} The isomorphism $T_X^*\Teich(S) \to T_X\Teich(S)$ determined by the Weil-Petersson metric is given by $\phi \mapsto \frac{\bar\phi}{\rho^2_X}$. Therefore the second statement follows by Theorem \ref{Rvolvariation}.

To see that $V$ is defined for all time we observe that in the Teichm\"uller metric, the norm of $V$ is bounded by $3/2$ from the Kraus-Nehari bound (Theorem \ref{kn_bound}). Since the Teichm\"uller metric is complete, a bounded vector field has a flow for all time.\eproof

Recall that $\mathcal V_R(N)$ is the infimum of the renormalized volume of $M \in CC(N)$. We define $\mathcal V_C(N)$ to be the same quantity with renormalized volume replaced by convex core volume. 
While $\mathcal V_C(N)$ is trivially non-negative, this is not clear for $\mathcal V_R(N)$. As noted above, if $N$ is a handlebody then work of Bridgeman-Canary \cite{BCrenorm} implies that  $\mathcal V_R(N) = -\infty$. However if $N$ has incompressible boundary we prove the following:

\begin{theorem}\label{renorm_inf}
Let $N$ be compact hyperbolizable 3-manifold with non-empty incompressible boundary and without torus boundary components. Then $\mathcal V_R(N) = \mathcal V_C(N)$. If there exists an $M \in CC(N)$ with $\Rvol(M) = \mathcal V_R(N)$ then either $N$ is acylindrical  and $M$ is the unique manifold in $CC(N)$ whose convex core boundary is totally geodesic or $N$ has the homotopy of a closed surface and $\mathcal V_R(N) = M$ if and only if $M$ is a Fuchsian manifold.
\end{theorem}

{\bf Proof:} We first observe that, by the upper bound on renormalized
volume from Theorem \ref{L2_squeeze}, $\mathcal V_R(N)\le \mathcal
V_C(N)$.  If we have $M \in CC(N)$ with $\Rvol(M) = \mathcal V_R(N)$
then $M$ is critical point of $\Rvol$ and therefore by the variational formula (Theorem \ref{Rvolvariation}), $\phi_M = 0$. This
occurs exactly when the convex core of $M$ has totally geodesic boundary which implies
that either $N$ is acylindrical or $M$ is Fuchsian. In the acylindrical case there is a unique $M \in CC(N)$ whose convex core boundary is totally geodesic.

Now choose $M \in CC(N)$ and let $M_t$ be the flow of $V$ with $M= M_0$ and let $\phi_t$ be the Schwarzian derivative of the projective boundary $M_t$.
We have
$$\Rvol(M_T) = \Rvol(M) - \int_0^T \|\phi_t\|^2_2 dt.$$
Since $\Rvol$ is bounded below on $CC(N)$ the integral $\int^\infty_0 \|\phi_t\|_2^2 dt$ converges. Therefore  there is a increasing sequence $\{t_i\}$ such that $t_i \rightarrow \infty$ and $\|\phi_{t_i}\|_2 \to 0$ as $i\to\infty$. We also have that $\Rvol(M_t)$ is a decreasing function of $t$ that is bounded below and hence $\Rvol(M_t)$ is convergent  as $t\to\infty$. Together with Theorem \ref{L2_squeeze} these two facts imply that
$$\underset{i\to\infty}{\lim} \Rvol(M_{t_i}) = \underset{i\to\infty}{\lim}\Cvol(M_{t_i}).$$
Since $\Rvol(M_{t_i})$ is a decreasing sequence we have
$$\Rvol(M) \ge \underset{i\to\infty}{\lim} \Rvol(M_{t_i}).$$
By definition $\Cvol(M_t) \ge \mathcal V_C(N)$ so
$$\underset{i\to\infty}{\lim}\Cvol(M_{t_i}) \ge \mathcal{V}_C(N).$$
Therefore $\Rvol(M) \ge \mathcal V_C(N)$. Since $M$ is arbitrary we have $\mathcal V_R(N) \ge \mathcal V_C(N)$ completing the proof. \eproof

By a theorem of Storm, (\cite[Theorem 5.9]{storm2}), the infimum of the volume of the convex core is half the simplicial volume of the double of the manifold with the infimum realized if and only $N$ is acylindrical or $N$ has the homotopy type of a closed surface. As an immediate corollary of our result and Storm's theorem we have:

\medskip

\noindent
{\bf Corollary \ref{simplicial_vol}} {\em
Let $N$ be a compact hyperbolizable 3-manifold with non-empty incompressible boundary and without torus boundary components. Then $\mathcal V_R(N) = \frac12 V_S(DN)$ where $DN$ is the double of $N$ and $V_S(DN)$ is the simplicial volume. The infimum is realized if and only if $N$ is acylindrical or has the homotopy type of a closed surface.}

\medskip
The manifold $DN$ is hyperbolic if and only if $N$ is acylindrical and then $V_S(DN)$ is twice the volume of the convex core of the unique $M\in CC(N)$ with totally geodesic boundary. As noted in the introduction, Palette has proved Corollary \ref{simplicial_vol} if $N$ is acylindrical. Palette's proof does not appeal to Storm's result so combining Theorem \ref{renorm_inf} together with Palette's work gives a new proof of the Storm theorem in the acylindrical case. In fact, by studying the limit of the $M_t$ as $t\to\infty$ one could directly prove Storm's theorem without appealing to \cite{PAL}. This will be discussed further in \cite{BBB}.

\bibliography{bib}
\bibliographystyle{math}
\end{document}